%% file: Physical_Realizability_Conditions_for_Mixed_Bilinear-Linear_Quantum_Cascades_with_Pure_Field_Coupling.tex
\documentclass[10pt,letterpaper,conference]{ieeeconf}

\input{CSSPreamble.tex}

\usepackage{cite} 
\usepackage{color} 
\usepackage{graphicx} 
\usepackage{latexsym,amssymb,amsmath} 
\usepackage{mathrsfs} 
\usepackage{multicol}
\usepackage{subeqn}
\usepackage[german,english]{babel}
\usepackage[T1]{fontenc}
\usepackage[latin1]{inputenc} 
\usepackage{widetext} 

\allowdisplaybreaks[4]

\IEEEoverridecommandlockouts                              
\title{Physical Realizability Conditions for Mixed Bilinear-Linear Quantum Cascades with Pure Field Coupling$^\ast$\thanks{This work was supported under
Australian Research Council's Discovery Projects funding scheme (project number DP110102322).}}

\author{Luis~A.~Duffaut~Espinosa$^\dagger$\thanks{$^\dagger$School of Engineering and Information Technology, University of New South Wales at ADFA, Canberra,
ACT 2600, Australia. {\tt\small \{l.duffaut,  i.petersen, v.ougrinovski\}@adfa.edu.au.}}, Z.~Miao$^\ddag$\thanks{$^\ddag$Research School of Information
Sciences and Engineering, Canberra, ACT 2601, Australia. {\tt\small zibo.miao@anu.edu.au}.}, I.~R.~Petersen$^\dagger$, V.~Ugrinovskii$\,^\dagger$, and
M.~R.~James$^\S$\thanks{$^\S$ARC Centre for Quantum Computation and Communication Technology, Research School of Engineering, Australian National University,
Canberra, ACT 0200, Australia. {\tt\small matthew.james@anu.edu.au}.}}

\begin{document}

\maketitle

\begin{abstract}
This paper aims to provide conditions under which a quantum stochastic differential equation can serve as a model for interconnection of a bilinear system 
evolving on an operator group $\pmb{SU(2)}$ and a linear quantum system representing a quantum harmonic oscillator. To answer this question we derive algebraic 
conditions for the preservation of canonical commutation relations (CCRs) of quantum stochastic differential equations (QSDE) having a subset of system 
variables satisfying the harmonic oscillator CCRs, and the remaining variables obeying the CCRs of $\pmb{SU(2)}$. Then, it is shown that from the physical 
realizability point of view such QSDEs correspond to bilinear-linear quantum cascades. 
\end{abstract}

\section{Introduction}  \label{sec:section1}

\begin{center}
\begin{figure}[h]
\begin{center}
\includegraphics[width=9cm]{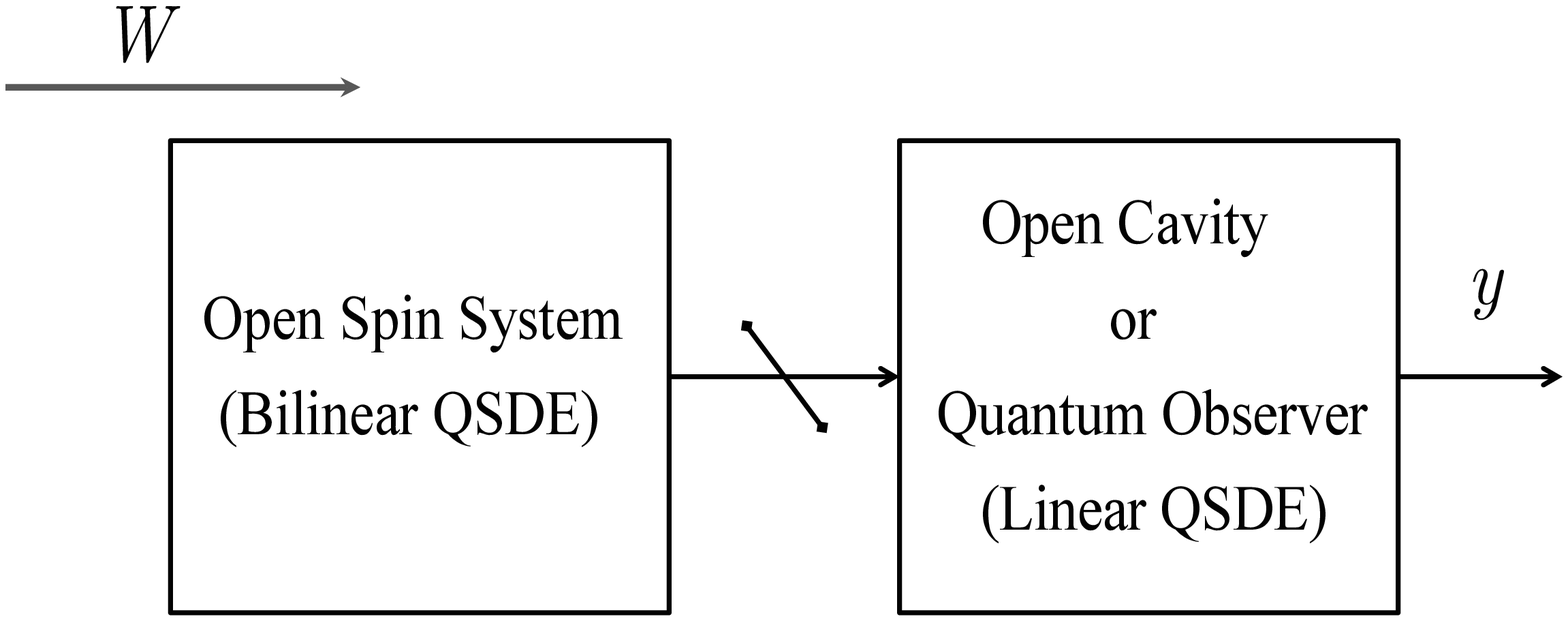}
\caption{Non interacting bilinear-linear quantum cascade open to a field $W$.}
\label{fig:BilinearLinearCascade}
\end{center}
\vspace{-0.3in}
\end{figure}
\end{center}

In many applications, systems are interconnected in order to form more complex systems. Open quantum systems are not the exception. For instance,
non-classical propagating electromagnetic fields, as now experimentally realizable, are an important resource in linear optics quantum information processing
\cite{Bozyigit_et_al_2011}. They can be constructed by cascading a two-level quantum system, as a source, with a cavity (quantum harmonic operator system) 
which filters the signals from the two-level system. In this case, the two-level system and the oscillator are separated by a transmission line such that there 
is no direct interaction between their system variables \cite{daSilva_et_all_2010} (Figure \ref{fig:BilinearLinearCascade}). 
From a control perspective, such apparatus are of great importance. For instance, a natural question is whether it is possible to estimate the states of a 
source system via a simpler oscillator system, the latter playing a role of a Luenberger observer. The answer to such question is by no means obvious, and it 
primarily depends on how one choses to describe the quantum nature of the comprising systems and the interconnection itself.   

It has been established that the framework of QSDEs provides an alternative description for studying quantum systems, in which it
allows the translation of standard control techniques into a quantum mechanical framework 
\cite{James-Nurdin-Petersen_2008,Dong-Petersen_2010,Belavkin_83,Helon-James_2006,Lloyd_2000,Maalouf-Petersen_2011,Sarovar-Ahn-Jacobs-Milburn_2004,
Yanagisawa-Kimura_2003, Yanagisawa-Kimura_2003a,Wiseman-Milburn_2010}. The QSDE description is in agreement with the \emph{Heisenberg picture} of 
quantum systems \cite{Parthasarathy_92}. Not every QSDE describes a quantum system (for instance, CCRs are not satisfied necessarily), however there exist 
conditions under which linear and bilinear QSDEs obey quantum mechanical laws, namely \emph{physical realizability conditions} 
\cite{James-Nurdin-Petersen_2008,Duffaut-et-al_2012b,Duffaut-et-al_2013a}. Physical realizability conditions provide simple testable matrix conditions 
containing the essentials for a system to be considered quantum. In this context, quantum oscillators are described by linear QSDEs and two-level 
systems are described by bilinear QSDEs. However, the the task of, for example, observing a physically realizable two level system with a physically realizable 
linear QSDE by cascading requires first of all to ensure the physical realizability of the composite system. Such cascade system goes beyond the realm in which 
the physical realizability of linear and bilinear QSDEs has been studied so far. Therefore, it is important to consider \emph{mixed} physical realizability 
conditions. That is to say, it is required a testable condition for the physical realizability of cascade bilinear-linear systems having a subset of system 
variables satisfying the harmonic oscillator CCRs, and the remaining variables obeying the CCRs of a two level system (i.e., the CCRs of ${SU(2)}$ 
\cite{Duffaut-et-al_2012b,Mahler-Weberrus_98}). An analysis of this type also provides a glimpse of the full characterization of bilinear QSDEs with additive 
and multiplicative quantum noise as open quantum systems.  

The earliest work on a systematic approach to cascade quantum systems can be trace to \cite{Gardiner_93,Carmichael_93}. In \cite{Gough-James_2009}, the 
treatment of the quantum cascading problem was extended in a manner that completely characterizes the dynamics of the composite system from a network point of 
view. This setting is natural from the engineering point of view where the decomposition of systems plays a fundamental role in systems analysis and synthesis. 
This approach has been proved valuable since it shows explicitly the interacting field channels, and hence interconnections via those channels can be 
constructed in a natural manner. In contrast, the more standard way of describing quantum systems via evolution of a density operator does not allow a network 
methodology explicitly, because the interacting channels are averaged out and therefore the interconnection cannot be described directly. One way to keep track 
of the information about the coupling channels is through the Belavkin filter \cite{Belavkin_83}, but this approach requires measurements such as homodyne or 
heterodyne detection~\cite{Wiseman-Milburn_2010}. Using such measurements is precluded when the objective is coherent control, i.e., when the controller or 
observer is itself a quantum system \cite{Lloyd_2000}. Still the approach in \cite{Gough-James_2009} starts from a purely quantum description to then using 
QSDEs to give the description of the cascade in terms of quantum operators, which is the opposite to what physical realizability conditions provide. In other 
words, it is desired for control applications to find conditions under which a cascaded QSDE preserves the physical realizability conditions of the composite 
systems (quantum coherent cascades, in our case), and therefore allow to identify the underlying quantum operators, when they exist, governing the dynamics of 
the cascade. In this regard, the goal of this paper is twofold. First, the aim is to obtain conditions for the preservation of physical realizability of 
bilinear QSDEs having both additive and multiplicative quantum noise inputs, and having initial conditions satisfying mixed CCRs (a combination between the 
harmonic oscillator and finite level systems CCRs). The second goal is to provide necessary and sufficient conditions for the physical realizability of the 
bilinear-linear cascade of QSDEs.

The paper is organized as follows. Section \ref{sec:2} presents the basic preliminaries on open quantum systems, in particular, harmonic oscillator systems, two
level systems and cascade of systems. In Section \ref{sec:3}, the algebraic machinery is given. This is followed by Section \ref{sec:4}, in which the result on
the preservation of mixed CCRs for bilinear QSDE with additive and multiplicative noise is developed. In Section \ref{sec:5}, the physical realizability of
bilinear-linear QSDE cascades is analyzed. Finally, Section \ref{sec:conclusions} gives the conclusions and future research directions to follow.

\section{Open Quantum Systems and their Cascade} \label{sec:2}

\subsection{Notation}

Let $\re$ denote the real numbers and $\C$ the complex numbers with imaginary unit $\pmb{i}$. The set of real and complex $n$-dimensional vectors are denoted 
$\re^n$ and $\C^n$, respectively. The set of real and complex $n$ by $m$ matrices are denoted $\re^{n\times m}$ and $\C^{n\times m}$. The $n$-dimensional 
identity matrix is denoted by $I_n$, and the $n \times m$ dimensional zero matrix is $0_{n\times m}$. A separable Hilbert space 
is denoted by $\mathfrak{H}$. The set of operators in $\mathfrak{H}$ is denoted by $\mathfrak{T}(\mathfrak{H})$, the set of $n$ dimensional vectors of 
operators 
in $\mathfrak{T}(\mathfrak{H})$ is denoted by $\mathfrak{T}(\mathfrak{H})^{n}$ and the set of $n\times m$ dimensional arrays of operators in 
$\mathfrak{T}(\mathfrak{H})$ is denoted by $\mathfrak{T}(\mathfrak{H})^{n\times m}$. The operator $\hat{I}$ denotes the identity in 
$\mathfrak{T}(\mathfrak{H})$. The operation $[\cdot,\cdot]:\mathfrak{T}(\mathfrak{H}) \times\mathfrak{T}(\mathfrak{H}) \rightarrow \mathfrak{T}(\mathfrak{H})$ 
is known as the {commutator}, and it is defined as $[x,y] = xy-yx$. For vectors $x\in \mathfrak{T}(\mathfrak{H})^n$ and $y\in 
\mathfrak{T}(\mathfrak{H})^m$ the commutator is given as
\begdi
[{x},y^T] \triangleq {x} y^T - (y x^T)^T \in \mathfrak{T}(\mathfrak{H})^{n\times m}, 
\enddi 
$ {x}^\# \triangleq \left( x_1^\ast \; x_2^\ast \; \hdots \; x_n^\ast \right)^T$, $x^\dagger=(x^\#)^T$, $(\cdot)^T$ denotes the transpose operation and 
$(\cdot)^\ast$ denotes the adjoint (or the complex conjugate
in the case of complex vectors or matrices). On a quantum mechanical framework, it is common to multiply either vectors or matrices by arrays of operators. For 
example, let $A\in \C^{m \times n}$ and ${X}\in \mathfrak{T}(\mathfrak{H})^{n \times m}$, 
the $(i,j)$ element of the multiplication of a matrix by an operator matrix is 
\begdi
(A X)_{ij} = \sum_{k=1}^n a_{ik} x_{kj}  \in \mathfrak{T}(\mathfrak{H}). 
\enddi
obeys the usual matrix multiplication rules. These considerations allow to treat operators as system variables since in quantum mechanics they play the role 
of states, and therefeore allow us to use state space systems notation.\\

\begre
The operations between complex matrices and operators follow the guidelines of the standard \emph{canonical quantization} [5], which in simple 
words is a recipe that promotes the system variables from a classical mechanical framework into an operator framework in order to obtain a quantum mechanical 
description of the system.
\endre

\subsection{Open quantum systems}

Quantum systems interacting with an external environment are known as \emph{open quantum systems}. \emph{Observables} in a Hilbert
space $\mathfrak{H}$ represent physical quantities that can be measured, while quantum states give the current status of the system. Here open quantum systems 
are treated in the context of quantum stochastic processes \cite{Bouten-Handel-James_2007,Parthasarathy_92}. The non-commutativity of observables is a 
fundamental difference between quantum systems and classical systems in which the former must satisfy certain CCRs, which lead to the \emph{Heisenberg 
uncertainty principle} \cite{Kennard_27}. The environment consists of a collection of oscillator systems, each with the annihilation field operator $w(t)$ and 
the creation field operator $w^\ast (t)$ used for annihilation and creation of quanta at point $t$, and commonly known as the \emph{boson quantum field} (a 
quantum version of a Wiener process). Here it is assumed that $t$ is a real time parameter. These operators generate three interacting signals in the evolution 
of the system: the annihilation processes $W(t)$, the creation process $W^\dagger (t)$, and 
the counting process $\Lambda(t)$. 

The unitary evolution of an observable $X \in \mathfrak{T}(\mathfrak{H})$ in the \emph{Heisenberg picture} is described by the operator equation
\begeq \label{eq:heisenberg_picture}
X(t)=U^\dagger(t) (X\otimes \hat{I})\, U(t), 
\endeq 
where $U(t)$ is unitary for all $t$, and is the solution of the operator stochastic differential equation
\begin{align*}
dU(t)  = &\left(\rule{0in}{0.18in} (S - \hat{I})\,d\Lambda(t) + L \,d W^\dagger(t) - L^\dagger S\,d W(t)\right. \\ 
& {}\;\;\; \left. - \frac{1}{2}( L^\dagger L + \pmb{i} \mathcal{H})\,dt \right) U(t),
\end{align*}
with initial condition $U(0) = \hat{I}$. $\mathcal{H}$ denotes the system \emph{Hamiltonian} of the system, and $L$ and $S$ (unitary) determine the 
\emph{coupling} of the system to the field and the interaction between fields, respectively. For simplicity, this paper will consider only one interactiong 
field $W$. Using the \emph{quantum It\^o formula} for $X_1,X_2 \in \mathfrak{T}(\mathfrak{H})$ \cite{Hudson_parthasarathy_84}, i.e.
\begeq \label{eq:Ito_formula}
d(X_1X_2) = (dX_1)X_2+X_1(dX_2) +(dX_1)(dX_2),
\endeq 
the dynamics of \rref{eq:heisenberg_picture} is expressed as
\begeq \label{eq:general_evolution}
\begin{split}
dX  = &   \, (S^\dagger X S-X)\,d \Lambda + {\cal L}(X) \,dt  + S^\dagger[X,L] \,dW^\dagger \\
& {} + [L^\dagger,X]S \, dW,
\end{split}
\endeq
where ${\cal L}(X)$ is the Lindblad operator defined as
\begeq \label{eq:Linblad_operator}
{\cal L}(X)=-{\pmb i}[X,\mathcal{H}] +\frac{1}{2}\left(L^\dagger[X,L] + [L^\dagger,X]L\right). 
\endeq 
The output field is given by $Y(t)=U(t)^{\dagger }W(t)U(t)$, which amount to 
\begeq \label{eq:general_output}
dY = L dt + S dW.
\endeq 
The dynamics of an open quantum systems is usually parametrized by the triple $(S,L,\mathcal{H})$. Henceforth assume that $S=\hat{I}$. 

It is often convenient to express QSDEs in terms of quadrature fields, which make all system matrices real. This is provided by the following linear 
transformation of the interacting fields 
\begeq \label{eq:quadrature_transformation}
\left(\begin{array}{c}
\bar{W}_1 \\ \bar{W}_2
\end{array}\right)= \left(\begin{array}{cc}
1 & 1 \\ -\pmb{i} & \pmb{i} 
\end{array}\right) 
\left(\begin{array}{c}
W \\ W^\dagger
\end{array}\right),
\endeq
where the operators $\bar{W}_1$ and $\bar{W}_2$ are now self-adjoint. Moreover, the It\^o table (see \cite{Hudson_parthasarathy_84}) for these quadrature fields
is
\begeq \label{eq:ito_table_quadrature}
\left(\begin{array}{c}
d\bar{W}_1 \\ d\bar{W}_2
\end{array}\right)    \left(\begin{array}{cc}
d\bar{W}_1 & d\bar{W}_2
\end{array}\right)= \left(\begin{array}{cc}
1 & \pmb{i} \\ -\pmb{i} & 1 
\end{array}\right) dt.
\endeq
Similarly, the quadrature form of the output fields can be obtained from the same quadrature transformation. Thus, 
\vspace*{0.05in}\begeq \label{eq:bilinear_system_output}
\colvec{2}{d Y_1}{d Y_2}= \colvec{2}{L + L^\#}{\pmb{i}(L^\# - L)} \,dt + \colvec{2}{d\bar{W}_1}{d\bar{W}_2}. 
\endeq

\subsection{Linear open quantum systems}

The Hilbert space for this class of systems is
$\mathfrak{H}_1=\ell^2(\C)$ (the space of square integrable complex sequences) \cite{Dong-Petersen_2010}, and the vector of 
system variables is 
\begeq \label{eq:x1}
x_1 \in \mathfrak{T}(\mathfrak{H}_1)^{2n},
\endeq 
For instance, a single harmonic oscillator system variables in terms of the annihilator operator $a$ and creation operator $a^\dagger$ is written in
self-adjoint form $x_1\in  \mathfrak{T}(\mathfrak{H}_1)^{2}$ by using the transformation
\begeq \label{eq:selfadjoint_aadagger}
x_1 = \left(\begin{array}{cc} 1 & 1 \\ -\pmb{i} & \pmb{i} \end{array} \right) \colvec{2}{a}{a^\dagger}.
\endeq
The CCRs for $a$ and $a^\dagger$ are $[a,a]=[a^\dagger,a^\dagger]=0$ and $[a,a^\dagger]=1$. For a vector of $n$ creation and $n$ annihilator operators, 
one has that
\begdi
\left[x_1, x_1^T\right] = \left[\colvec{5}{a_1}{a_1^\dagger}{\vdots}{a_n}{a_n^\dagger}, (a_1\;a_1^\dagger \; \hdots \;a_n \;a_n^\dagger)\right] = (I_{n} 
\otimes J),
\enddi
where 
\begdi
J = \left( \begin{array}{cc} 0 & 1 \\ - 1 & 0 \end{array} \right).
\enddi
In self-adjoint form, by applying \rref{eq:selfadjoint_aadagger}, the CCRs are
\begeq \label{eq:CCRs_x1}
[x_1,x_1^T]= 2\pmb{i} \underbrace{(I_{n} \otimes J)}_{\displaystyle \triangleq \Theta}.
\endeq
The Hamiltonian for this class of systems is the quadratic form $\mathcal{H}_1= x_1^T R x_1$ with $R$ real symmetric, and the coupling operator is considered to
be linear, i.e., $L_1 = \Gamma_1 x_1$. The general form for the QSDE having these Hamiltonian an coupling operator is
\begin{subequations} \label{eq:linear_QSDE}
\begin{align}
\label{eq:linear_QSDE_evolution} dx_1   & = A x_1\, dt + B\, d\bar{W}  \\
\label{eq:linear_QSDE_output} dy_1  & = C x_1\, dt +  d\bar{W},
\end{align}
\end{subequations}
where $A \in \re^{3\times 3}$, $B\in \re^{n\times 2}$ and $C \in \re^{2\times n} $, and $\bar{W} = (\bar{W}_1 \;\bar{W}_2 )^T$.

For system \rref{eq:linear_QSDE} to have any hope of being quantum mechanical, it is fundamental that system \rref{eq:linear_QSDE} preserves \rref{eq:CCRs_x1} 
over time. The next theorem gives conditions for the preservation of CCRs of $x_1$ over time. 
\begth (See \cite{James-Nurdin-Petersen_2008,Dong-Petersen_2010}.) 
\label{th:preservingCCRx1}
QSDE \rref{eq:linear_QSDE_evolution} with system variables as in \rref{eq:x1} satisfying $[x_1(0),x_1(0)^T] = 2\pmb{i}\Theta$ implies $[x_1(t),x_1(t)^T] =
2\pmb{i}\Theta$ for all $t\ge 0$ if and only if
\begeq \label{eq:linear_CCR_preservation}
A\Theta + \Theta A^T + BJB^T =0,
\endeq
\endth

\subsection{Two level open quantum system}

For an open two-level quantum system interacting with one boson quantum field, the Hilbert space is $\mathfrak{H}_2=\C^2$ and the vector of system variables is 
\begeq \label{eq:x2}
x_2 \in \mathfrak{T}(\mathfrak{H}_2)^3,
\endeq
Note that operators in $\mathfrak{T}(\mathfrak{H}_2)$ are simply matrices in $\C^{2\times 2}$. These operators are chosen to be self-adjoint, so that $x_2$ 
satisfies $x_2=x_2^\#$.~In particular, an operator $\hat{\sigma} \in \mathfrak{T}(\mathfrak{H}_2)$ is spanned by the Pauli matrices \cite{Mahler-Weberrus_98}, 
i.e., 
$\hat{\sigma}=\frac{1}{2}\sum_{i=0}^3 \kappa_i \sigma_i$, 
where $\kappa_0=\Tr(\hat{\sigma})$, $\kappa_i=\Tr(\hat{\sigma}\sigma_i)$, and
\begdi
\sigma_0=\left(\begin{array}{cc}
         1 & 0 \\ 0 & 1
        \end{array} \right),\;\; 
\sigma_1=\left(\begin{array}{cc}
         0 & 1 \\ 1 & 0
        \end{array} \right),\\
\enddi
\begdi
\sigma_2=\left(\begin{array}{cc}
         0 & -{\pmb i} \\ {\pmb i} & 0
        \end{array} \right),\;\;
\sigma_3=\left(\begin{array}{cc}
         1 & 0 \\ 0 & -1
        \end{array} \right)
\enddi
denote the Pauli matrices. Thus, $\kappa_0,\kappa_1,\kappa_2$ and $\kappa_3$ determine uniquely the operator $\hat{\sigma}$. The product of Pauli matrices
satisfy 
\begeq \label{eq:Pauli_CCR_multiplication}
\sigma_i\sigma_j = \delta_{ij}I_3+ {\pmb i} \sum_{k}\epsilon_{ijk}\sigma_k,
\endeq 
and therefore its CCRs are
\begeq \label{eq:Pauli_CCR} 
[\sigma_i,\sigma_j] = 2{\pmb i} \sum_{k}\epsilon_{ijk}\sigma_k,
\endeq
where $\delta_{ij}$ is the Kronecker delta and $\epsilon_{ijk}$ denotes the Levi-Civita tensor. Given that \rref{eq:Pauli_CCR_multiplication} allows to 
write any product Pauli operators as linear forms, a large class of polynomial quantum systems can be characterized by considering linear 
Hamiltonian and 
coupling operators, i.e., 
${\mathcal{H}_2}=\alpha_2 x_2$ and 
$L_2=\Gamma_2 x_2$, 
where $\alpha_2^T\in \re^3$ and $\Gamma_2^T \in \C^3$.

Observe that, in general, the evolution of $x_2$ is a bilinear QSDEs with only multiplicative quantum noise expressed as 
\vspace*{-0.1in}\begin{subequations} \label{eq:bilinear_QSDE}
\begin{align} 
\label{eq:bilinear_system}      & dx_2 =A_0\,dt+Ax_2 \, dt +{B}_1x_2\,d\bar{W}_1 + {B}_2x_2\, d\bar{W}_2, \\
\label{eq:bilinear_QSDE_output} & dy_2 = C x_2\, dt + d\bar{W},
\end{align}
\end{subequations}
where $A_0\in \re^3$, $A,{B}_1,{B}_2\in \re^{3\times 3}$ and $C \in \re^{2\times n} $. Conditions for CCR preservation of $x_2$ are given in the next theorem.

\begth (See \cite{Duffaut-et-al_2012b,Duffaut-et-al_2013a}.) \label{th:preservingCCRx2} QSDE \rref{eq:bilinear_system} with system
variables as in \rref{eq:x2} satisfying $[x_2(0),x_2(0)^T] = 2\pmb{i}\Theta^-(x_2(0))$ implies $[x_2(t),x_2(t)^T] = 2\pmb{i}\Theta^-(x_2)$ for all $t\ge 0$ if
and only if
\begin{subequations}
\begin{align}
\label{eq:CCRx21} {B}_1 + {B}_1^T = {B}_2 + {B}_2^T & = 0 \\
\label{eq:CCRx22} {B}_1{B}_2^T - {B}_2{B}_1^T -\Theta(A_0)  & = 0 \\
\label{eq:CCRx23} A^T+A + {B}_1 {{B}_1}^T + {B}_2 {{B}_2}^T & = 0.
\end{align} 
\end{subequations}
\endth

The fact that all matrices in systems \rref{eq:linear_QSDE} and \rref{eq:bilinear_QSDE} are real is due to the quadrature
transformation \rref{eq:quadrature_transformation}.
\subsection{Cascades of open quantum systems}

If the cascade connection of a two level system and a linear quantum system is considered, the composite system lives in $\mathfrak{H}_{12} = \mathfrak{H}_{1}
\otimes \mathfrak{H}_{2} = \ell^2(\C) \otimes \C^2$, which is the completion of the direct product of $\ell^2(\C)$ and $\C^2$. In this construction the system
variables in $\mathfrak{H}_1$ and $\mathfrak{H}_2$ when embedded in $\mathfrak{H}_{12}$ commute between each other. The cascade of open quantum systems is
described by an algebraic operation on the $(S,L,\mathcal{H})$ parametrization. Such operation is defined next. 

\begde (See \cite{Gough-James_2009}.) Given two open quantum systems parametrized by ${\bf G_1}=(S_1,L_1,\mathcal{H}_1)$ and ${\bf G_2}=(S_2,L_2,\mathcal{H}_2)$
having the same number of field channels, the \emph{series product}  ${\bf G_1} \triangleleft {\bf G_2}$ is defined as 
\begin{align} \label{eq:composite_SLH}
\begin{split}
{\bf G_1} \triangleleft {\bf G_2} = & {\,} \left(\rule{0in}{0.2in}S_2S_1,L_2+L_1, \right. \\
                                    &  {\;}\left.\mathcal{H}_1 + \mathcal{H}_2 + \frac{1}{2\pmb{i}}\left( L_2^\dagger S_2 L_1 - L_1^\dagger S_2^\dagger L_2
\right)\right)
\end{split}
\end{align}
\endde
Since we assume both $S_1$ and $S_2$ to be the identity operators in the corresponding spaces, the QSDE describing the cascade of systems (system $2$
drives system $1$) can then be written for $x^T=(x_1^T\; x_2^T)$ as 
\begin{align} \label{eq:cascade_QSDE}
\begin{split}
\hspace{-0.11in}dx \!= & \left({\cal L}_{1}(x)+{\cal L}_{2}(x) + L_2^\dagger[x,L_1] + [L_1^\dagger,x]L_2\right)dt \\
                                        &  + [x,L_2+L_1]\,dW + [L_2^\dagger+L_1^\dagger,x]\,dW^\dagger.
\end{split}
\end{align}

\section{Some Algebraic Relations} \label{sec:3}

Let $\beta=(\beta_1,\beta_2,\beta_3)^T \in \C^3$, and define the linear mapping $\Theta^-: \C^3 \rightarrow \C^{3\times 3}$ such that
\begdi
\Theta^-(\beta)= \left(\begin{array}{ccc}
         0 & \beta_3 & -\beta_2 \\ -\beta_3 & 0 & \beta_1 \\ \beta_2  & -\beta_1  & 0
        \end{array} \right).
\enddi  
This mapping is understood for vector of operators by associating with $\beta$ the vector of operators $\hat{\beta}=(\beta_1 \hat{I},\beta_2 \hat{I},\beta_3
 \hat{I})^T\in \mathfrak{T}(\mathfrak{H}_2)^3$ such that 
\begdi
\Theta^-(\hat{\beta})= \left(\begin{array}{ccc}
         0 & \beta_3 \hat{I} & -\beta_2\hat{I} \\ -\beta_3\hat{I} & 0 & \beta_1\hat{I} \\ \beta_2\hat{I}  & -\beta_1\hat{I}  & 0
        \end{array} \right) \in \mathfrak{T}(\mathfrak{H}_2)^{3\times 3},
\enddi 
Abusing the notation, $\hat{I}$ will be omitted hereafter, and the fact that $\beta$ is either a vector of numbers or a vectors of operators will be 
understood from the context. As an example, the product of Pauli operators can be expressed in a compact matrix form thanks to the mapping $\Theta^-(\cdot)$. 
That is, 
\begdi 
x_2x_2^T=I_{3} + {\pmb i} \Theta^-(x_2) \in \mathfrak{T}(\mathfrak{H}_2)^{3\times 3}.
\enddi 
Observe here that the identity matrix $I_3$, under our convention, is strictly speaking denoting a three dimensional diagonal matrix of the identity operator 
in $\mathfrak{T}(\mathfrak{H}_2)$. Similarly, the CCRs for Pauli operators are written as 
\begdi
[x_2,x_2^T]=2 {\pmb i} \Theta^-(x_2) \in \mathfrak{T}(\mathfrak{H}_2)^{3\times 3}. 
\enddi
Considering the \emph{stacking operator}, denoted $\vec$, whose action on an $m\times n$ dimensional array creates a $mn$ dimensional column vector by stacking 
its columns below one another. Applying $\vec$ to $\Theta^-(\beta)$ gives 
$\vec(\Theta^-(\beta)) = F \beta$, 
where $m=n=3$, $F\triangleq \left(  F_1 , F_2 , F_3  \right)^T$, the $(j,k)$ component of $F_i$ is $(F_i)_{jk}=\epsilon_{ijk}$, and $\epsilon_{ijk}$ is the
Levi-Civita tensor. Some properties of $\Theta^-(\cdot)$ are summarized in the next lemma (see \cite{Duffaut-et-al_2013a} for more identities). 
\begle (See \cite{Duffaut-et-al_2012b,Duffaut-et-al_2013a}.) \label{le:theta_properties}
The mapping $\Theta^-(\cdot)$ satisfies 
\renewcommand*\theenumi{{\roman{enumi}}}
\renewcommand*\labelenumi{$\theenumi.$}
\begin{enumerate}
\item \label{eq:Theta_1} $\Theta^-(\beta)\gamma = - \Theta^-(\gamma) \beta,$ \vspace*{0.05in}

\item \label{eq:Theta_beta_beta} $\Theta^-(\beta)\beta = 0,$ \vspace*{0.05in}

\item \label{eq:Theta_composition} $\Theta^-\left(\Theta^-(\beta)\gamma\right)=[\Theta^-(\beta),\Theta^-(\gamma)]$.
\end{enumerate}
\endle \vspace*{0.05in}
This properties hold when $\beta$ and $\gamma$ are either $\C^3$ vectors or $\mathfrak{T}(\mathfrak{H}_2)^{3}$ vectors.

The explicit computation of the vector fields in \rref{eq:general_evolution} and \rref{eq:cascade_QSDE} for $x_1$ and $x_2$ is given in the next lemma. 
\begle \label{le:vector_fields} The nonzero coefficients of equations \rref{eq:general_evolution} and \rref{eq:Linblad_operator} for the dynamics of $x_1$, 
$x_2$ and
the cascade $G_1 \triangleleft G_2$ are 
\begin{align*}
   [x_1,\mathcal{H}_1]  & = 2\pmb{i} \Theta R x_1, \\ 
  [x_1,L_1]  & = 2\pmb{i} \Theta \Gamma_1^T, \\  
  [x_1,L_1^\dagger]  & = 2\pmb{i} \Theta \Gamma_1^\dagger, \\  
  L_1^\dagger[x_1,L_1]  & = 2\pmb{i} \Theta \Gamma_1^T\Gamma^\# x_1,\\
  [x_1,L_1^\dagger]L_1  & = -2\pmb{i} \Theta \Gamma_1^\dagger\Gamma x_1,   \\
  [x_2,\mathcal{H}_2]  & = -2{\pmb i}\Theta^-(\alpha_2^T)x_2, \\
  [x_2,L_2]  & = -2{\pmb i}\Theta^-(\Gamma_2^T)x_2, \\
  [x_2,L_2^\dagger]  & = -2{\pmb i}\Theta^-(\Gamma_2^{\dagger})x_2, \\
  L_2^\dagger[x_2,L_2]  & = -2\pmb{i}\Theta^-(\Gamma_2^T)\Gamma_2^\dagger+2\Theta^-(\Gamma_2^T)\Theta^-(\Gamma_2^\dagger)x_2, \\
 [x_2,L_2^\dagger]L_2   & = 2{\pmb i}\Theta^-(\Gamma_2^T)\Gamma_2^\dagger-2\Theta^-(\Gamma_2^\dagger)\Theta^-(\Gamma_2^T)x_2, \\
 L_2^\dagger [x_1,L_1] & = 2\pmb{i} \Theta \Gamma_1^T \Gamma_2^\# x_2, \\
 [L_1^\dagger,x_1]L_2  & = -2\pmb{i} \Theta \Gamma_1^\dagger \Gamma_2  x_2.
\end{align*}
\endle \vspace*{0.1in}

From this lemma, system \rref{eq:linear_QSDE} is written as
\begeq \label{eq:physical_linear_evolution}
\begin{split}
dx_1  = {} & 2 \Theta\left(R+ \mathfrak{F}(\Gamma_1^\dagger \Gamma_1) \right) x_1\,dt  \\
      & \hspace*{-0.2in} + 2\pmb{i}\Theta \left( \left(-\Gamma_1^\dagger + \Gamma_1^T\right) \;\; -\pmb{i}\left(\Gamma_1^\dagger + \Gamma_1^T\right)
\right)\,d\bar{W}, \\
dy_1  = {} & \colvec{2}{\Gamma_1 + \Gamma_1^\#}{i(\Gamma_1^\# - \Gamma_1)} x_1 \,dt + d \bar{W},  
\end{split}
\endeq
where $\mathfrak{F}(z) \triangleq \frac{1}{2\pmb{i}}\left(z-z^* \right)$ is the imaginary part of $z$. \\

\begre We see from  (\ref{eq:physical_linear_evolution}) that a linear coupling operator $L_1$ produces, in $\mathcal{L}_1(x_1)$, only linear terms of the form
$M x_1\, dt$ with $M \in \C^{2n\times 2n}$, and constant noise vector fields because of the CCRs of $x_1$. Suppose now that $L_1$ is a quadratic form, i.e.,
$L_i=x_1^T \Gamma_1 x_1$, then the term $[L_1^\dagger,x_1]$ produces a bilinear term, however evaluating, for instance, $[L_1,x_1]L_1^\dagger$ generates a term
of the form $M_1 (x_1\otimes x_1)$ with $M_1 \in \C^{2n\times (2n)^2}$. Even more, these terms cannot be embedded in a higher dimensional bilinear system since
by doing so only produces polynomials of higher order of the oscillator system variables. This indicates that a QSDE describing a system of $n$ harmonic
oscillators cannot have terms of the form $B_ix_1 d\bar{W}_i$ when the coupling operator is a linear form. 
\endre

\noindent For system \rref{eq:bilinear_QSDE}, one has that
\begeq \label{eq:physical_spin_evolution}
\begin{split}
\hspace*{-0.1in} dx_2  = {} &  -2{\pmb i}\Theta^-(\Gamma_2^T)\Gamma_2^\dagger \,dt -2 \Theta^-(\alpha_2^T)x_2\,dt \\
      & \hspace*{-0.2in} +\left(\Theta^-(\Gamma_2^T)\Theta^-(\Gamma_2^\dagger) +\Theta^-(\Gamma_2^\dagger)\Theta^-(\Gamma_2^T)\right)x_2\,dt \\
      & \hspace*{-0.2in} +\pmb{i} \Theta^-(\Gamma_2^\dagger - \Gamma_2^T) x \,d\bar{W}_1 - \Theta^-(\Gamma_2^T + \Gamma_2^\dagger) x_2 \,d\bar{W}_2, \\
dy_2  = {} & \colvec{2}{\Gamma_2 + \Gamma_2^\#}{i(\Gamma_2^\# - \Gamma_2)} x_2 \,dt + d \bar{W}.     
\end{split}
\endeq
Finally, \rref{eq:cascade_QSDE} for the cascade of \rref{eq:bilinear_QSDE} driving \rref{eq:linear_QSDE} is
\begeq \label{eq:physical_cascade_evolution}
\begin{split}
\hspace*{-0.1in}\colvec{2}{dx_1}{dx_2} & = \colvec{2}{0}{-2{\pmb i}\Theta^-(\Gamma_2^T)\Gamma_2^\dagger }dt \\
     & \hspace*{-0.25in} + \left(\begin{array}{cc} \mathcal{R}_1    &   -4\Theta\mathfrak{F}(\Gamma_1^T\Gamma_2^\#) \\ 0  & \mathcal{R}_2 \end{array}  
\right)\colvec{2}{x_1}{x_2} \, dt \\
     & \hspace*{-0.25in} + \left( \begin{array}{cc}
             0 & 0 \\
             0 & \pmb{i} \Theta^-(\Gamma_2^\dagger - \Gamma_2^T) \end{array} \right) \colvec{2}{x_1}{x_2} \, d\bar{W}_1 \\
     & \hspace*{-0.25in} - \left(\begin{array}{cc}
             0 & 0 \\
             0 & \Theta^-(\Gamma_2^T + \Gamma_2^\dagger) \end{array} \right) \colvec{2}{x_1}{x_2} \, d\bar{W}_2 \\
     & \hspace*{-0.25in} +  \left(\begin{array}{c}
             2\pmb{i}\Theta \left( \left(-\Gamma_1^\dagger + \Gamma_1^T\right) \;\; -\pmb{i}\left(\Gamma_1^\dagger + \Gamma_1^T\right)\right)  \\  0 \end{array}
\right) d\bar{W}
\end{split}
\endeq
with $\mathcal{R}_1 = 2 \Theta\left(R+ \mathfrak{F}(\Gamma_1^\dagger \Gamma_1) \right)$ and $\mathcal{R}_2 = -2 
\Theta^-(\alpha_2^T) + \Theta^-(\Gamma_2^T)\Theta^-(\Gamma_2^\dagger) +\Theta^-(\Gamma_2^\dagger)\Theta^-(\Gamma_2^T)$.

We observe that the QSDE (\ref{eq:physical_cascade_evolution}) contains both additive and multiplicative noise terms, and its drift term is affine. Two
question can now be asked. The first is under what conditions a general QSDE of such form (see equation (\ref{eq:bilinearQSDE_additive_multiplicative_system})
below) preserves the CCRs for $x_1$ and $x_2$ at the same time. This question is addressed in Section \ref{sec:4}. Then, it will be desired to know under what
conditions there exists $(S,L,\mathcal{H})$ as in \rref{eq:composite_SLH} such that \rref{eq:cascade_QSDE} can be written as in
\rref{eq:general_evolution} (Section~\ref{sec:5}). 

\section{Preservation of CCRs} \label{sec:4}

Consider an arbitrary $n$-dimensional bilinear QSDE interacting with a quadrature field. That is,
\begin{subequations} \label{eq:bilinearQSDE_additive_multiplicative_system}
\begin{align}  
\label{eq:bilinearQSDE_additive_multiplicative} & \hspace*{-0.15in} dx = A_0dt + Ax dt + B_1 x d\bar{W}_1+B_2 x d\bar{W}_2 + B d\bar{W},\\
\label{eq:bilinearQSDE_additive_multiplicative_output} & \hspace*{-0.15in} dy = C xdt + d\bar{W},
\end{align}
\end{subequations}
where $A_0 \in \re^n$, $A,B_1,B_2 \in \re^{n\times n}$, $B \triangleq(\bar{B}_1\; \bar{B}_2)$, $\bar{B}_1, \bar{B}_2 \in \re^n$, and $d\bar{W}= (d\bar{W}_1 
\; d\bar{W}_2)^T$. 

In previous work (\cite{James-Nurdin-Petersen_2008,Duffaut-et-al_2013a}), the quantum noise appearing in the equations was either additive or multiplicative.
This model differs from those in what it includes both additive and multiplicative noise, and the system models are such that their system variables can be
partitioned into two mutually commuting sets each having different CCRs. Specifically, one set obeys the CCRs of harmonic oscillators, and the other follows the
CCRs of a two-level system. That is,
\begin{align} \label{eq:mixed_CCRs}
\nonumber [x,x^T] = & \left[\colvec{2}{x_1}{x_2}, (x_1^T \; x_2^T) \right] \\
         = & \left( \begin{array}{cc} \Theta & 0 \\ 0 & \Theta^-(x_2) \end{array} \right).
\end{align}
Conversely, the imposition of these CCRs on an arbitrary $x$ induces automatically a partition of $x$ in a way that one set obeys harmonic oscillator CCRs,
while the other obey the CCRs of $SU(2)$. Since this partition of $x$ can always be obtained via a linear transformation, one can assume without loss of 
generality that $x$ is always of the form $x^T=(x_1^T \; x_2^T)$. 

Consider now the block partition of $A_0$, $A$, $B_i$ and $\bar{B}_i$ as follows
\begin{align*}
& A_0 = \colvec{2}{A_{01}}{A_{02}}, \;\; A = \left(\begin{array}{cc} A_{11} &  A_{12} \\ A_{21}  & A_{22} \end{array}  \right),\\ 
& B_i = \left(\begin{array}{cc} B_{i11} &  B_{i12} \\ B_{i21}  &  B_{i22} \end{array}  \right) \;\;{\rm and}\;\; \bar{B}_i  =
\colvec{2}{\bar{B}_{i1}}{\bar{B}_{i2}}
\end{align*}
for $i=1,2$. Recalling the fact that $x_2$ is self-adjoint, one can infer that 
\begin{align*}
\bar{B}_{i2} = 0_{3 \times 1}. 
\end{align*}
This agrees with the fact that a bilinear QSDE is driving a linear QSDE. In summary, the only source of additive noise is provided by the linear QSDE.
Note that the bilinear QSDE system can only provide multiplicative noise to the composite system. Also, the equation for $dx_1$ can only have bilinear terms
with respect to $x_2$. This means that 
\begin{align*}
B_{i12} = 0_{2\times2}. 
\end{align*}

\begth \label{th:preservationCCR_linear-bilinear_cascade} Let $x$ be a vector of operators satisfying CCRs \rref{eq:mixed_CCRs}, a QSDE as in
\rref{eq:bilinearQSDE_additive_multiplicative} preserves such CCRs for all $t \ge 0$ if and only if the linear QSDE 
\begdi dx_1 = A_{11} x_1\, dt + ( \bar{B}_{11}\; \bar{B}_{21}) \, d\bar{W} \enddi and the bilinear QSDE 
\begdi dx_2 = (A_{02} + A_{22}x_2)\,dt + B_{122}x_2 \, d\bar{W}_1 + B_{222}x_2 \, d \bar{W}_2 \enddi
satisfy the conditions in Theorems \ref{th:preservingCCRx1} and \ref{th:preservingCCRx2}, respectively, in addition to 
\begin{align} \label{eq:consistency_condition1}
\nonumber \lefteqn{\hspace*{-0.4in}(I_{3}\otimes A_{12})F + (B_{122}^T\otimes \bar{B}_{21}) }\\
& \hspace*{1in} - (B_{222}^T \otimes \bar{B}_{11})  =0. 
\end{align}
\endth \vspace*{0.1in}

\begre The structure showed in \rref{eq:physical_cascade_evolution} appears naturally from the preservation of mixed CCRs (see the proof of Theorem 
\ref{th:preservationCCR_linear-bilinear_cascade} in the appendix). 
\endre

\section{Cascade Physical Realizability}   \label{sec:5}

As mentioned in the introduction, physical realizability for linear and bilinear QSDEs has previously been treated independently of each other
(\cite{James-Nurdin-Petersen_2008,Duffaut-et-al_2012b,Duffaut-et-al_2013a}). However, a more natural setting for quantum systems is when linear and $n$-level
systems are components of a larger system. The objective here is to give conditions for physical realizability for a bilinear QSDE driving a linear
QSDE. The general notion of physical realizability is provided next. It basically ties QSDE's of arbitrary nature with an $(S,L,\mathcal{H})$ parametrization.
\begde 
\label{def:physical realizability} A QSDE is said to be \emph{physically realizable} if there exist operators $\mathcal{H}$ and $L$ such that the QSDE can be
written as in \rref{eq:general_evolution} and \rref{eq:general_output}.
\endde

In what follows a summary of the necessary and sufficient conditions for linear and bilinear QSDE's is given. Then the second main result of the paper is
given. That is,  necessary and sufficient conditions for physical realizability of the cascade of a bilinear QSDE followed by a linear QSDE.  

\subsection{Physical realizability of linear QSDEs}

\begde \label{def:physical realizability_x1} The system \rref{eq:linear_QSDE} is said to be physically realizable if there exist $\mathcal{H}_1$ and $L_1$ such
that \rref{eq:linear_QSDE} can be written as in \rref{eq:general_evolution} and \rref{eq:general_output}.
\endde
The explicit form of matrices $A, B, C_1$ and $C_2$ in \rref{eq:linear_QSDE} is given in terms of a Hamiltonian and coupling operator next, and can be
identified from \rref{eq:physical_linear_evolution}. The existence of an $(S_1,L_1,\mathcal{H}_1)$ parametrization of linear QSDEs with system variables as in
\rref{eq:x2} is given by the next theorem.
\begth (See \cite{James-Nurdin-Petersen_2008,Dong-Petersen_2010}.) \label{th:physical_realizability_linear}
System \rref{eq:linear_QSDE} is physically realizable if and only if 
\vspace*{0.05in}
\renewcommand*\theenumi{{\roman{enumi}}}
\renewcommand*\labelenumi{$\theenumi.$}
\begin{enumerate}
\item \label{itm:linear_PR1} $A\Theta + \Theta A + B J B  = 0$,
\vspace*{0.05in}
\item \label{itm:linear_PR2} $B  =  \Theta C^T (J \otimes I_{n})$,
\vspace*{0.05in}
\end{enumerate}
where $\mathcal{H}_1$ and $\Gamma_1$ are uniquely identified as
\begdi
R = \frac{1}{4}\left( -\Theta A + A^T \Theta \right) \;\; {\rm and}\;\; \Gamma_1 = \frac{1}{2} \left(C_1+\pmb{i} C_2 \right). 
\enddi
\endth

Note that $\rref{itm:linear_PR1}$ is identical to \rref{eq:linear_CCR_preservation}, however the latter is generated purely form algebraic 
considerations. 

\subsection{Physical realizability of bilinear QSDEs}

\begde \label{def:physical realizability_x2} System \rref{eq:bilinear_QSDE} is said to be
{physically realizable} if there exist $\mathcal{H}$ and $L$ such that \rref{eq:bilinear_system} can be written as in \rref{eq:general_evolution} and
\rref{eq:general_output}.
\endde

The explicit matrices $A_0,A, B_{1}, B_{2},C_1$ and $C_2$ in terms of a Hamiltonian and coupling operator can be extracted from
\rref{eq:physical_spin_evolution}. The existence of an $(S_2,L_2,\mathcal{H}_2)$ parametrization of bilinear QSDEs with system variables as in \rref{eq:x2} is
given by the next theorem.

\begth (See \cite{Duffaut-et-al_2012b,Duffaut-et-al_2013a}.) \label{th:physical_realizability_bilinear}
The system \rref{eq:bilinear_QSDE} with output equation \rref{eq:bilinear_system_output} is physically realizable if and only if 
\vspace*{0.05in}
\renewcommand*\theenumi{{\roman{enumi}}}
\renewcommand*\labelenumi{$\theenumi.$}
\begin{enumerate}
\item \label{itm:PR_bilinear1} $\displaystyle A_0=\frac{1}{2}({B}_{1} + \pmb{i}{B}_{2}) \left(C_1+\pmb{i} C_2 \right)^\dagger$,
\vspace*{0.05in}
\item \label{itm:PR_bilinear2} $\displaystyle {B}_{1}= \Theta^-(C_2^T)$,
\vspace*{0.1in}
\item \label{itm:PR_bilinear3} $\displaystyle {B}_{2}= -\Theta^-(C_1^T)$,
\vspace*{0.05in}
\item \label{itm:PR_bilinear4} $\displaystyle A+A^T+  B_1 B_1^T + B_2 B_2^T= 0$.
\vspace*{0.05in}
\end{enumerate}
In which case, one can identify the matrix $\alpha_2$ defining the system Hamiltonian and the coupling matrix $\Gamma_2$ as
\begdi
  \alpha_2= \frac{1}{8}\vec(A-A^T)^T  F, \;\; {\rm and} \;\; \Gamma_2=\frac{1}{2}(C_1+\pmb{i}C_2).
\enddi
\endth 

Similar to the case of linear QSDEs, condition $\rref{itm:PR_bilinear4}$ is identical to \rref{eq:CCRx23}, however \rref{eq:CCRx23} is obtained form purely 
algebraic considerations.

\subsection{Physical realizability of a class of cascade bilinear-linear QSDE's}

The second main result of the paper is now presented. First, the definition of a physically realizable bilinear-linear cascade is given. 
\begde \label{de:cascade_PR}
\label{def:physical realizability_cascade} A QSDE is said to be a \emph{physically realizable bilinear-linear cascade} if there exist operators $\mathcal{H}$
and $L$ as in \rref{eq:composite_SLH} such that QSDE \rref{eq:cascade_QSDE} can be written as in \rref{eq:general_evolution} and \rref{eq:general_output}.
\endde

The characterization of the physical realizability of a bilinear-linear cascade of QSDEs is given in the next theorem. 
\begth \label{th:cascade_linear-bilinear_PR} The system \rref{eq:bilinearQSDE_additive_multiplicative_system} is physically realizable according to Definition
\ref{de:cascade_PR} if and only if the following conditions hold 
\renewcommand*\theenumi{{\roman{enumi}}}
\renewcommand*\labelenumi{$\theenumi.$}
\begin{enumerate}
\item \label{itm:PR_cascade1} The matrices $A_0$, $A$, $B_1$, $B_2$, $B$ and $C$ in \rref{eq:bilinearQSDE_additive_multiplicative_system} are of the following 
form
 \begin{align*}
& \hspace*{-0.3in} A_0 = \colvec{2}{0}{A_{02}},\;\;\; A  = \left(\begin{array}{cc} A_{11} &  A_{12} \\ 0  & A_{22} \end{array}  \right), \\ 
& \hspace*{-0.3in} B_i  = \left(\begin{array}{cc} 0 &  0 \\ 0  &  B_{i22} \end{array}  \right), \;\;\;
 B  = \left(\begin{array}{cc}  \bar{B}_{11} & \bar{B}_{21} \\  0 & 0 \end{array} \right), \\
& \hspace*{-0.3in}{\rm and}\;\; C= \colvec{2}{C_{1}}{0}.
\end{align*}

\item \label{itm:PR_cascade2} System 
\begin{align*} dx_1 & = A_{11} x_1\, dt + ( \bar{B}_{11}\; \bar{B}_{21}) \, d\bar{W},\\
dy & = C_1 x_1 + d\bar{W} \end{align*} 
is physically realizable in the sense of Definitions~\ref{def:physical realizability_x1}. 

\item \label{itm:PR_cascade3} System 
\begin{align*} dx_2 & = (A_{02} + A_{22}x_2)\,dt\! + \! B_{122}x_2 \, d\bar{W}_1 \! +\! B_{222}x_2 \, d \bar{W}_2,\\
dy & = C_2 x_2 + d\bar{W} \end{align*} 
is physically realizable in the sense of Definitions~\ref{def:physical realizability_x2}, where $C_2^T = (C_{21}^T \; C_{22}^T)$ is such that the following 
consistency condition holds: 
\begin{align} 
\label{eq:composition_consistency}
A_{12} & = \bar{B}_{11} C_{21}+\bar{B}_{21} C_{22}.
\end{align}
\end{enumerate}
\endth \vspace*{0.1in}

The following corollary is a consequence of the previous theorem.
\begco \label{co:CCR_imply_PR} A bilinear-linear cascade physically realizable QSDE preserves \rref{eq:mixed_CCRs}. 
\endco

\section{Conclusions and Future Research} \label{sec:conclusions}

Conditions for the preservation of mixed CCRs were developed. In particular, these conditions were obtained for bilinear systems having both additive and
multiplicative quantum noise inputs. It was also shown that bilinear-linear QSDE cascades are physical realizable when the linear and bilinear
subsystems are physically realizable and a consistency condition holds.

A future research direction is to consider an interactive Hamiltonian in the formalism (a hermitian operator $\mathcal{H}_I = x_1^T R_1 x_2$). This would 
allow our theory to capture some of the commonly used models in quantum optics. For example, an atom trapped in an optical cavity is described by the 
Jaynes-Cummings model, i.e. a model with a Hamiltonian of the form 
\begdi
\mathcal{H} = \frac{\hbar}{2} \omega_0 \sigma_z + \frac{1}{2} \gamma a \sigma^+ + \frac{1}{2} \gamma^* a^\dagger \sigma^- + \hbar \omega_c a^\dagger a, 
\enddi
where $\omega_c$ and $\omega_0$ are the frequencies of the cavity and atom, respectively, and $\gamma$ is the interaction strength. In addition, the 
conditions provided in this manuscript will potentially allow the synthesis of coherent quantum observers for $n$-level systems in the Heisenberg picture.

\section*{Acknowledgement}

The authors want to thank M. Wooley for useful discussions and insight on the physics relevance of the results presented in this paper.

\appendix[Proofs of Results]

\emph{Proof of Theorem \ref{th:preservationCCR_linear-bilinear_cascade}:} Using 
\rref{eq:Ito_formula} and \rref{eq:ito_table_quadrature}, it follows that $d[x,x^T]$ can 
be obtained by computing $d(xx^T)$ and $(d(xx^T))^T$. That is, 
\begin{align*}
\lefteqn{\hspace*{-0.13in}d(xx^T) = (dx)x^T+x(dx)^T+(dx)(dx)^T}\\
	= \,& {} (A_0 x^T + x A_0^T)\,dt + (Axx^T+xx^TA^T)\,dt \\ 
	  & {}  + ({B}_1xx^T+xx^T{B}_1^T)\,d \bar{W}_1 +({B}_2xx^T+xx^T{B}_2^T)\,d\bar{W}_2  \\
	  & {}  + {B}_1x x^T {B}_1^T d\bar{W}_1 d\bar{W}_1+{B}_1x x^T {B}_2^T d\bar{W}_1 d\bar{W}_2 \\
	  & {}  + {B}_2x x^T {B}_1^T d\bar{W}_2 d\bar{W}_1 +{B}_2x x^T {B}_2^T d\bar{W}_2 d\bar{W}_2\\
	  & {}  + (\bar{B}_1  d\bar{W}_1 + \bar{B}_2 d\bar{W}_2) (B_1 x  d\bar{W}_1 + B_2 x d\bar{W}_2)^T \\
  	  & {}  + (B_1 x  d\bar{W}_1 + B_2 x d\bar{W}_2)(\bar{B}_1  d\bar{W}_1 + \bar{B}_2 d\bar{W}_2)^T \\ 
  	  & {}  + (\bar{B}_1  d\bar{W}_1 + \bar{B}_2 d\bar{W}_2)(\bar{B}_1  d\bar{W}_1 + \bar{B}_2 d\bar{W}_2)^T \\
	  = \,& {} (A_0 x^T + x A_0^T)\,dt + (Axx^T+xx^TA^T)\,dt \\ 
	  & {}  + ({B}_1xx^T+xx^T{B}_1^T)\,d \bar{W}_1 +({B}_2xx^T+xx^T{B}_2^T)\,d\bar{W}_2  \\
	  & {}  + {B}_1x x^T {B}_1^T dt + \pmb{i}{B}_1x x^T {B}_2^T dt \\
	  & {}  - \pmb{i}{B}_2x x^T {B}_1^T dt + {B}_2x x^T {B}_2^T dt\\
	  & {}  + \bar{B}_1 x^T B_1^T\, dt + \pmb{i} \bar{B}_1 x^T B_2^T\, dt\\
	  & {}  - \pmb{i}\bar{B}_2 x^T B_1^T\, dt + \bar{B}_2 x^T B_2^T\, dt\\
	  & {}  +  B_1 x \bar{B}_1^T\, dt + \pmb{i} B_1 x \bar{B}_2^T\, dt\\
	  & {}  - \pmb{i} B_2 x \bar{B}_1^T\, dt + B_2 x \bar{B}_2^T\, dt\\
	  & {}  + \bar{B}_1 \bar{B}_1^T\, dt - \pmb{i} \bar{B}_1 \bar{B}_2^T\, dt\\
	  & {}  + \pmb{i} \bar{B}_2 \bar{B}_1^T\, dt + \bar{B}_2 \bar{B}_2^T\, dt.
	  \end{align*}
Similarly, 
\begin{align*}
\lefteqn{\hspace*{-0.18in}\left(d(xx^T)\right)^T}\\
	= \,& {} (A_0 x^T + x A_0^T)\,dt + (A(xx^T)^T + (xx^T)^TA^T)\,dt \\ 
	  & {}  + ({B}_1(xx^T)^T + (xx^T)^T{B}_1^T)\,d \bar{W}_1 \\
	  & {}  +({B}_2(xx^T)^T + (xx^T)^T{B}_2^T)\,d\bar{W}_2  \\
	  & {}  + {B}_1 (xx^T)^T {B}_1^T dt + \pmb{i}{B}_1 (xx^T)^T {B}_2^T dt \\
	  & {}  - \pmb{i}{B}_2 (xx^T)^T {B}_1^T dt + {B}_2 (xx^T)^T {B}_2^T dt\\
	  & {}  + B_1 x \bar{B}_1^T \, dt + \pmb{i}   B_2 x \bar{B}_1^T\, dt\\
	  & {}  - \pmb{i} B_1 x \bar{B}_2^T\, dt +  B_2 x \bar{B}_2^T\, dt\\
	  & {}  +  \bar{B}_1 x^T B_1^T\, dt + \pmb{i} \bar{B}_2 x^TB_1^T\, dt\\
	  & {}  - \pmb{i} \bar{B}_1 x^T B_2^T\, dt +  \bar{B}_2 x^T B_2^T\, dt\\
	  & {}  + \bar{B}_1 \bar{B}_1^T\, dt - \pmb{i}  \bar{B}_2 \bar{B}_1^T \, dt\\
	  & {}  + \pmb{i}  \bar{B}_1 \bar{B}_2^T\, dt + \bar{B}_2 \bar{B}_2^T \, dt.
\end{align*}
Hence, the commutator dynamics is
\begin{align*}
d\left[x,x^T\right] = & {\;} A[x,x^T] + [x,x^T] A^T     \\ 
	    & {} + (B_1[x,x^T] + [x,x^T] B_1^T ) d\bar{W}_1 \\
    	    & {} + (B_2[x,x^T] + [x,x^T] B_2^T ) d\bar{W}_2 \\
    	    & {} + (B_1[x,x^T]B_1^T + B_2[x,x^T] B_2^T ) \,dt \\
    	    & {} + \pmb{i}(B_1\{x,x^T\}B_2^T - B_2\{x,x^T\} B_1^T ) \,dt \\
            & {} + 2\pmb{i} B_2 x \bar{B}_1^T \, dt + 2\pmb{i} B_1 x  \bar{B}_2^T \, dt \\
            & {} + 2\pmb{i} \bar{B}_1 x^T B_2^T \, dt + 2\pmb{i} \bar{B}_2 x^TB_1^T \, dt \\
            & {} + \pmb{i} \bar{B}_1 \bar{B}_2^T \,dt - \pmb{i}\bar{B}_2 \bar{B}_1^T \, dt.
\end{align*}
To preserve \rref{eq:mixed_CCRs}, \rref{eq:bilinearQSDE_additive_multiplicative} has to satisfy 
\begeq \label{eq:commutaor_dynamic}
d\left[x,x^T\right] = 2 \pmb{i}\left(\begin{array}{cc} 0 & 0 \\
0 & \Theta^-(d x_2) \end{array}\right), 
\endeq 
where 
\begin{align*}
\Theta^-(d x_2) = & \; \Theta^-\left(A_{02}\right)\,dt + \Theta^-\left(A_{22} x_2\right)\, dt \\
& \, + \left(\Theta^-\left(B_{121} x_1 \right)+\Theta^-\left(B_{122}x_2 \right)\right)\, d\bar{W}_1\\
& \, + \left(\Theta^-\left(B_{221}x_1 \right)+\Theta^-\left( B_{222}x_2\right)\right)\, d\bar{W}_2.
\end{align*}
$A_{01}$ does not play a role in the preservation of CCRs. Therefore, without loss of generality $A_{01}$ is assumed to be zero. This goes
in agreement with the fact that no term of this type is generated by quantum systems originating from harmonic oscillators of the class considered in this
paper.

From \cite[Proposition 27.3]{Parthasarathy_92}, one can also equate the integrands in \rref{eq:commutaor_dynamic} to zero. Recall that $x_2(0)$ is represented
by the complete orthonormal set. This implies that any linear combination $\sum_{k=0}^s a_i x_i(0) \not = 0$ unless $a_i=0$ for all $i$ and $a_i \in \C$. In
addition, no linear combination of Pauli matrices generates $I_{3}$. Therefore, any equation  $Ax_2=b$ ($A\in \C^{3\times 3}$ and $b\in \C^3$) implies
$A=0$ and
$b=0$. These facts are summarized in the following equations that have to be satisfied for the preservation of CCRs.

\begin{widetext}
\begin{subequations}
\begin{align}
\label{subeq:preservation0} &B_{122} B_{222}^T - B_{222} B_{122}^T - \Theta^-(A_{02}) =0,  \\
\label{subeq:preservation1} &B_{i21} \Theta =0,  \\
\label{subeq:preservation2} &B_{i12} \Theta^-(x_2) = 0\\
\label{subeq:preservation3} &B_{i22} \Theta^-(x_2) +\Theta^-(x_2) B_{i22}^T  - \Theta^-(B_{i22}x_2)  = 0, \\
\label{subeq:preservation4} &  A_{11} \Theta + \Theta A_{11}^T + \pmb{i} \left( \bar{B}_{11} \bar{B}_{21}^T - \bar{B}_{21} \bar{B}_{11}^T \right) =0 \\
\label{subeq:preservation5} &  A_{12}\Theta^-(x_2)  + \bar{B}_{21} x_2^T B_{122}^T - \bar{B}_{11} x_2^T B_{222}^T  =0 \\
\label{subeq:preservation6} &  A_{21} \Theta = 0 \\
\label{subeq:preservation7} &  A_{22}\Theta^-(x_2) + \Theta^-(x_2) A_{22}^T + B_{122}\Theta^-(x_2) B_{122}^T + B_{222}\Theta^-(x_2) B_{222}^T  -
\Theta^-(A_{22} x_2)=0.
\end{align}
\end{subequations}
\end{widetext}

Relations \rref{subeq:preservation0}, \rref{subeq:preservation3} and \rref{subeq:preservation7} provide the preservation of CCRs of $x_2$ (Theorem
\ref{th:preservingCCRx2}). Similarly, \rref{subeq:preservation4} assures the preservation of CCRs for $x_1$ (Theorem \ref{th:preservingCCRx1}). Relations
\rref{subeq:preservation1}, \rref{subeq:preservation2} and \rref{subeq:preservation6} impose a structure on the blocks in matrices $A$, $B_1$, $B_2$ and
$B$. That is, one has that \rref{subeq:preservation2} provides $B_{i12}=0$ by the linear independence on the components of $x_2$. Since $\Theta$ only
permutes the rows and columns of $A_{12}$ and multiplies some of its components by $-1$ then $A_{12}=0$. The same argument provides $B_{i21}=0$. From
Lemma $3$ in \cite{Duffaut-et-al_2013a}, \rref{subeq:preservation3} is always satisfied, and allows to write $B_{i22}=\Theta^-(b_i)$ with 
\begdi 
b_i = -\frac{1}{n} \colvec{3}{ \Tr(F_1 B_{i22}) }{ \vdots }{ \Tr(F_s B_{i22})}.
\enddi
Therefore, by fixing the CCRs of $x$, the matrices in \rref{eq:cascade_QSDE} assume naturally the following structure
\begin{align*}
&  A_0 = \colvec{2}{0}{A_{02}},\;\;\; A    = \left(\begin{array}{cc} A_{11} &  A_{12} \\ 0  & A_{22} \end{array}  \right), \\ 
&  B_i  = \left(\begin{array}{cc} 0 &  0 \\ 0  &  B_{i22} \end{array}  \right), \;\;{\rm and}\;\;
 B  = \left(\begin{array}{cc}  \bar{B}_{11} & \bar{B}_{21} \\  0 & 0 \end{array} \right).
\end{align*}
To obtain \rref{eq:consistency_condition1}, first recall $\vec(ABC) = (C^T \otimes A)\vec(B)$ for $A, B$ and $C$ of appropriate dimensions. Then, applying the
stacking operator to \rref{subeq:preservation5} the desired consistency condition \rref{eq:consistency_condition1} is obtained.

Conversely, since the steps used above to obtain \rref{eq:consistency_condition1} are reversible and the fact that the preservation of CCRs for $x_1$ and 
$x_2$
in Theorems \ref{th:preservingCCRx1} and \ref{th:preservingCCRx2} imply \rref{subeq:preservation0}, \rref{subeq:preservation3}, \rref{subeq:preservation4} and
\rref{subeq:preservation7}, then \rref{eq:commutaor_dynamic} holds. This finalizes the proof. 
\endpr

\emph{Proof of Theorem \ref{th:cascade_linear-bilinear_PR}:} If system \rref{eq:bilinearQSDE_additive_multiplicative_system} is bilinear-linear cascade 
physically realizable, then it can be written as in \rref{eq:physical_cascade_evolution}, and the systems formed by matrices $(A_{11},(\bar{B}_{11} \, 
\bar{B}_{21}), C_1)$ and $(A_{02},A_{22},B_{122},B_{222},C_2)$ can be written as in \rref{eq:physical_linear_evolution} and 
\rref{eq:physical_spin_evolution}, 
respectively. Therefore, $(I,L_1,\mathcal{H}_1)$ and $(I,L_2,\mathcal{H}_2)$ can be identified so that the parametrization $(S,L,\mathcal{H})$ as in 
\rref{eq:composite_SLH} holds. It is only left to prove that $A_{12}$ can be written as \rref{eq:composition_consistency}. One has from Lemma 
\ref{le:vector_fields} that
\begin{align*}
A_{12} & = -4\Theta \mathfrak{F}(\Gamma_1^T \Gamma_2^\#) = 2\pmb{i}\Theta\Gamma_1^T \Gamma_2^\# - 2\pmb{i}\Theta\Gamma_1^\dagger \Gamma_2,
\end{align*}
and that $2\pmb{i}\Theta\Gamma_1^T = \bar{B}_{11}+\pmb{i} \bar{B}_{21}$ and $-2\pmb{i}\Theta\Gamma_1^\dagger=\bar{B}_{11}+\pmb{i}\bar{B}_{21}$. Thus, 
\begin{align} \label{eq:A12_PR_proof}
\nonumber A_{12} & = (\bar{B}_{11}+\pmb{i}\bar{B}_{21}) \Gamma_2^\# + (\bar{B}_{11}+\pmb{i}\bar{B}_{21}) \Gamma_2 \\
\nonumber       & = \bar{B}_{11} (\Gamma_2+\Gamma_2^\#)+\bar{B}_{21} \pmb{i}(\Gamma_2^\#-\Gamma_2)\\
       & = \bar{B}_{11} C_{21}+\bar{B}_{21} C_{22}.
\end{align}

On the other hand, assuming $\rref{itm:PR_cascade1}$-$\rref{itm:PR_cascade3}$ hold, then from Theorems \ref{th:physical_realizability_linear} and 
\ref{th:physical_realizability_bilinear} the triples $(I,L_1,\mathcal{H}_1)$ and $(I,L_2,\mathcal{H}_2)$ are uniquely identified. In particular, $\Gamma_1 
=\frac{1}{2}( C_{11}+ \pmb{i}C_{12})$ and $\Gamma_2 = \frac{1}{2}(C_{21}+ \pmb{i}C_{22})$. Finally, since \rref{eq:composition_consistency} hold all the 
steps 
in \rref{eq:A12_PR_proof} are reversible, then $A_{12} = -4\Theta \mathfrak{F}(\Gamma_1^T \Gamma_2^\#)$ as in \rref{eq:physical_cascade_evolution}. This 
completes the proof.
\endpr

\emph{Proof of Corollary \ref{co:CCR_imply_PR}:} To prove this result using Theorem \ref{th:cascade_linear-bilinear_PR}, only condition 
\rref{eq:consistency_condition1} of Theorem~\ref{th:preservationCCR_linear-bilinear_cascade} need to be established. Given that the cascade is physically 
realizable, one has that $A_{12}= 2\pmb{i} \Theta ( \Gamma_1^T\Gamma_2^\# - \Gamma_1^\dagger \Gamma_2 )$, $\bar{B}_{11}=\pmb{i}\Theta 
(\Gamma_1^T-\Gamma_1^\dagger)$, $\bar{B}_{21} =  \Theta (\Gamma_1^T+\Gamma_1^\dagger)$, $B_{122}=\pmb{i}\Theta^-( \Gamma_2^\dagger - \Gamma_2^T)$ and 
$B_{222}=-\Theta^-( \Gamma_2^T +\Gamma_2^\dagger )$. Using Lemma \ref{le:theta_properties}, it then follows that 
\begin{align*}
\lefteqn{\bar{B}_{21} x_2^T B_{122}^T} \\
      & =  \pmb{i} \Theta( \Gamma_1^\dagger+\Gamma_1^T ) x_2^T \Theta^-( \Gamma_2^\dagger - \Gamma_2^T ) \\
      & = \pmb{i} \Theta \left( \Gamma_1^T \Gamma_2 - \Gamma_1^T \Gamma_2^\# + \Gamma_1^\dagger \Gamma_2
          -\Gamma_1^\dagger \Gamma_2^\# \right) \Theta^-(x_2),
\end{align*}
\begin{align*}
\lefteqn{\bar{B}_{11} x_2^T B_{222}^T} \\
      & = -  \pmb{i} \Theta( \Gamma_1^T-\Gamma_1^\dagger ) x_2^T \Theta^-( \Gamma_2^T +\Gamma_2^\dagger ) \\
      & = \pmb{i} \Theta \left( \Gamma_1^T \Gamma_2 + \Gamma_1^T \Gamma_2^\# - \Gamma_1^\dagger \Gamma_2
          -\Gamma_1^\dagger \Gamma_2^\# \right) \Theta^-(x_2).
\end{align*}
Hence
\begin{align*}
A_{12}\Theta^-(x_2) +  \bar{B}_{21} x_2^T B_{122}^T-\bar{B}_{11} x_2^T B_{222}^T =0, 
\end{align*}
which is equivalent to \rref{eq:consistency_condition1} after applying the stacking operator and using the linear independence of the components of $x_2$. 
\endpr

\end{document}

%% file: CSSPreamble.tex
\usepackage{latexsym, amssymb, amsmath}
\usepackage[dvips]{graphicx} 

\usepackage{float}      
\restylefloat{figure}

\usepackage{graphicx,dsfont}
\usepackage{mathrsfs}  
\usepackage{color} 
\usepackage[latin1]{inputenc}
\usepackage{enumerate}
\usepackage{flushend} 
\usepackage{stfloats}
\usepackage{cite} 

\usepackage{balance}

\allowdisplaybreaks[4] 

\newtheorem{prop}{Proposition} 
\newtheorem{proposition}[prop]{Proposition} 
\newtheorem{lem}{Lemma}
\newtheorem{lemma}[lem]{Lemma} 
\newtheorem{thm}{Theorem} 
\newtheorem{theorem}[thm]{Theorem} 
\newtheorem{cor}{Corollary} 
\newtheorem{corollary}[cor]{Corollary} 
\newtheorem{defn}{Definition} 
\newtheorem{definition}[defn]{Definition} 
\newtheorem{exmp}{Example} 
\newtheorem{example}[exmp]{Example} 
\newtheorem{exam*}{Example}

\def\custombibliography#1{
 \normalsize
\section*{\centering References}
 \list
 {[\arabic{enumi}]}{\settowidth\labelwidth{[#1]}\leftmargin\labelwidth
 \setlength{\itemsep}{.1em}
 \advance\leftmargin\labelsep
 \usecounter{enumi}}
 \def\newblock{\hskip .11em plus .33em minus -.07em}
 \sloppy
 \sfcode`\.=1000\relax}

\def\L2{{\cal L}_2}
\def\bull{\rule{0.08in}{0.08in}} 
\def\openbull{\framebox[0.08in][c]{$\;$}} 
\def\re{{\mathbb R}} 

\def\C{{\mathbb C}} 

\def\ss#1{{\scriptstyle #1}} 
\def\eqref#1{(\ref{#1})} 

\def\Tr{{\rm Tr}}
\def\vec{{\rm vec}}

\newcommand{\comment}[1]{} 

\def\begce{\begin{center}}
\def\endce{\end{center}}
\def\begar{\begin{array}}
\def\endar{\end{array}}
\def\begeq{\begin{equation}}
\def\endeq{\end{equation}}
\def\begdi{\begin{displaymath}}
\def\enddi{\end{displaymath}}
\def\begdis{\begin{eqnarray*}}
\def\enddis{\end{eqnarray*}}
\def\begeqa{\begin{eqnarray}}
\def\endeqa{\end{eqnarray}}
\def\begdes{\begin{description}}
\def\enddes{\end{description}}
\def\begit{\begin{itemize}}
\def\endit{\end{itemize}}
\def\begen{\begin{enumerate}}
\def\enden{\end{enumerate}}
\def\beglar{\left[\begin{array}}
\def\endrar{\end{array}\right]}
\def\begle{\begin{lemma}}
\def\endle{\end{lemma}}
\def\begde{\begin{definition}}
\def\endde{\end{definition}}
\def\begth{\begin{theorem}}
\def\endth{\end{theorem}}
\def\begco{\begin{corollary}}
\def\endco{\end{corollary}}
\def\begprop{\begin{proposition}}
\def\endprop{\end{proposition}}
\def\begex{\begin{example}}
\def\endex{\hfill\openbull \end{example} \vspace*{0.1in}}
\def\begexer{\begin{exercise}}
\def\endexer{\end{exercise}}
\def\begre{\noindent{\em Remark}: }
\def\endre{\\}
\def\begres{\noindent{\bf Remarks}:\begin{enumerate}}
\def\endres{\end{enumerate} \par}

\def\endpr{\hfill\bull \vspace*{0.1in}}

\def\begtab{\begin{tabular}}
\def\endtab{\end{tabular}}
\def\rref#1{(\ref{#1})}

\newcommand\cdcout[1]{} 

\newcommand{\rv}[1]{\boldsymbol{#1}} 
\newcommand{\RomanNumber}[1]{\uppercase\expandafter{\romannumeral #1}}
\newcommand{\romannumber}[1]{\lowercase\expandafter{\romannumeral #1}}

\DeclareMathAlphabet{\mathpzc}{OT1}{pzc}{m}{it}

\def\1{\rv 1} 

\def\allseriesA{\mbox{$\re\langle\langle A\rangle\rangle$}}
\def\allseriesA'{\mbox{$\re\langle\langle A'\rangle\rangle$}}

\def\L1spaceprodu{{ L}_1(\Omega\times [0,T],{\mathcal P},P\otimes \lambda)}

\def\Hspace0{{\mathcal H}^2_0}

\newcount\colveccount
\newcommand*\colvec[1]{
        \global\colveccount#1
        \begin{pmatrix}
        \colvecnext
}
\def\colvecnext#1{
        #1
        \global\advance\colveccount-1
        \ifnum\colveccount>0
                \\
                \expandafter\colvecnext
        \else
                \end{pmatrix}
        \fi
}





%% file: Physical_Realizability_Conditions_for_Mixed_Bilinear-Linear_Quantum_Cascades_with_Pure_Field_Coupling.bbl
\begin{thebibliography}{999} 


\bibitem{Belavkin_83} V.~Belavkin, ``On the theory of controlling observable quantum systems,'' \emph{Automation and Remote Control}, vol.~42, no.~2,
pp.~178--188, 1983.

\bibitem{Bouten-Handel-James_2007} L.~Bouten, R.~Van Handel, and M.~R.~James. ``An introduction to quantum filtering,'' \emph{SIAM Journal of Control and
Optimization}, vol.~46, no.~6, pp.~2199--2241, 2007.

\bibitem{Bozyigit_et_al_2011} D.~Bozyigit, C.~Lang, L.~Steffen, J.~M.~Fink, C.~Eichler, M.~Baur, R.~Bianchetti, P.~J.~Leek, S.~Filipp, M.~P.~da~Silva, A.~Blais,
and A.~Wallraff, ``Antibunching of microwave-frequency photons observed in correlation measurements using linear detectors,'' Nature Phys., vol.~7, pp.
154--158, 2011. 


\bibitem{Carmichael_93} H.~J.~Carmichael, ``Quantum trajectory theory for cascaded open systems,'' \emph{Phys. Rev. Lett.}, vol.~70, no.~15, pp.~2273--2276,
1993.

\bibitem{Dirac_25} P.~A.~M. Dirac, ``The Fundamental Equations of Quantum Mechanics". \emph{Proc. R. Soc. Lond. A}, vol.~109, pp.~642--653, 1925.

\bibitem{Helon-James_2006} C.~D'Helon and M.~R.~James, ``Stability, gain and robustness in quantum feedback networks,'' \emph{Physical Review
A}, vol.~73, pp.0533803, 2006

\bibitem{daSilva_et_all_2010} M.~Da Silva, D.~Bozyigit, A.~Wallraff, and A.~Blais, ``Schemes for the observation of photon correlation functions in circuit
QED with linear detectors.'' Phys. Rev. A, vol.~82, pp.~043804, 2010.

\bibitem{Doherty-Jacobs_99} A.~Doherty and K.~Jacobs, ``Feedback-control of quantum systems using continuous state-estimation,'' \emph{Physical Review
A}, vol.~60, pp. 2700--2711, 1999.

\bibitem{Dong-Petersen_2010} D.~Dong and I.~R.~Petersen, ``Quantum control theory and applications: A survey,'' \emph{IET Control Theory \&
Applications}, vol.~4, no.~12, pp.~2651--2671 2010.


\bibitem{Duffaut-et-al_2012b} L.~A.~Duffaut Espinosa, Z.~Miao, I.~R.~Petersen, V.~Ugrinovskii and M.~R.~James, ``Preservation of Commutation Relations and
Physical Realizability of Open Two-Level Quantum Systems,'' \emph{Proc. 51th IEEE Conference on Decision}, Maui, Hawaii, pp~3019--3023, 2012.

\bibitem{Duffaut-et-al_2013a} L.~A.~Duffaut Espinosa, Z.~Miao, I.~R.~Petersen, V.~Ugrinovskii and M.~R.~James, ``On the preservation of commutation and
anticommutation relations of $n$-level quantum systems,'' \emph{Proc. 2013 American Control Conference}, Washington D.C., 2013, to appear.

\bibitem{Gardiner_93} C.~W.~Gardiner, ``Driving a quantum system with the output field from another driven quantum system,'' \emph{Phys. Rev. Lett.},
vol.~70, no.~15, pp.~2269--2272, 1993.

\bibitem{Gough-James_2009} J,~Gough and M.~R.~James, ``The series product and its application to quantum feedforward and feedback networks,'' \emph{IEEE
Transactions on Automatic Control}, vol.~54, no.~11, pp.~2530--2544, 2009.

\bibitem{Hudson_parthasarathy_84} R.~L.~Hudson and K.~R.~Parthasarathy, ``Quantum It\^o Formula and Stochastic Evolutions,'' \emph{Communications in
Mathematical Physics}, vol.~93, pp.~301--323, 1984.


\bibitem{James-Nurdin-Petersen_2008}
M.~R.~James, H.~I.~Nurdin, and I.~R.~Petersen,  ``$H^\infty$ Control of linear quantum stochastic systems,'' \emph{IEEE Transactions on Automatic
Control}, vol.~53, pp.~1787--1803, 2008.

\bibitem{Kennard_27} 
E.~H.~Kennard, ``Zur Quantenmechanik einfacher Bewegungstypen'', \emph{Zeitschrift f\"ur Physik}, vol.~ 44, no.~4-5, pp.~326--352, 1927.


\bibitem{Lloyd_2000} S.~Lloyd, ``Coherent quantum feedback,'' \emph{Physical Review A}, vol.~62, pp.~022108, 2000.

\bibitem{Maalouf-Petersen_2011}
A.~I.~Maalouf and I.~R. Petersen, ``Bounded Real Properties for a Class of Annihilation-Operator Linear Quantum Systems,'' \emph{IEEE Transactions on Automatic 
Control}, 2011, vol.~56, pp.~786--801, 2011.

\bibitem{Mahler-Weberrus_98}
G.~Mahler and W. V. Weberu\ss, \emph{Quantum Networks: Dynamics of Open Nanostructures}, Springer, Berlin, 1998.

\bibitem{Parthasarathy_92}
K.~R.~Parthasarathy, \emph{An Introduction to Quantum Stochastic Calculus}, Birkh\"auser Verlag, Berlin, 1992.

\bibitem{Sarovar-Ahn-Jacobs-Milburn_2004} M.~Sarovar, C.~Ahn, C.~Jacobs, and G.~J.~Milburn, ``Practical scheme for error control using feedback,''
\emph{Physical Review A}, vol.~69, pp.052324, 2004. 

\bibitem{Wiseman-Milburn_2010}
H.~Wiseman and G.~Milburn, \emph{Quantum Measurement and Control}, Cambridge University Press, New York, 2010.


\bibitem{Yanagisawa-Kimura_2003} M.~Yanagisawa and H.~Kimura, ``Transfer function approach to quantum control-part I: Dynamics of quantum feedback
systems,'' \emph{IEEE Transactions on Automatic Control}, vol.~48, no.~12, pp. 2107--2120, 2003.

\bibitem{Yanagisawa-Kimura_2003a} M.~Yanagisawa and H.~Kimura, ``Transfer function approach to quantum control-part II: Control concepts and
applications,'' \emph{IEEE Transactions on Automatic Control}, vol.~48, no.~12, pp. 2121--2132, 2003.

\end{thebibliography}
